\documentclass[]{amsart}
\usepackage{amsmath,amssymb,graphicx}
\newtheorem{theo}{\sc Theorem}
\newtheorem{claim}{\sc Claim}
\newtheorem{lem}{\sc Lemma}
\newtheorem{rem}{\sc Remark}

\newtheorem{example}{\sc Example}
\def\C{\mathbb{C}}
\def\CD{\widehat\C}
\def\Z{\mathbb Z}

\def\ende{\square}
\def\Ende{$\ende$}
\def\be{\begin{equation}}
\def\ee{\end{equation}}
\def\ds{\displaystyle}
\def\ts{\textstyle}
\def\proof{\noindent {\it Proof.~}}
\def\und{\quad{\rm and}\quad}
\address{Institut f\"ur Mathematik $\cdot$ TU Dortmund $\cdot$ D-44221 Dortmund $\cdot$ Germany}
\email{stein@math.tu-dortmund.de}

\begin{document}
\begin{center} {\Large\sc Algebraic Curves and Meromorphic Functions\\ Sharing Pairs of Values\\}

\medskip \sc Norbert Steinmetz\end{center}

\begin{quote}{\sc Abstract.} The {\sc 4IM+1CM}-Problem is to determine all pairs $(f,g)$ of meromorphic functions in the complex plane that are not M\"obius transformations of each other and share five pairs of values, one of them {\sc CM} (counting multiplicities). In the present paper it is shown that each such pair parameterises some algebraic curve $K(x,y)=0$ of genus zero
and bounded degree. Thus the search may be restricted to the pairs of meromorphic functions $(Q(e^z),\widetilde Q(e^z))$, where $Q$ and $\widetilde Q$ are non-constant rational functions
of low degree. This leads to the paradoxical situation that the {\sc 4IM+1CM}-problem could be solved by a computer algebra virtuoso rather than a complex analyst.\end{quote}


\section{Introduction}

Transcendental meromorphic functions $f$ and $g$ are said to share the {\it pair} $(a,b)$ of complex numbers if $f(z)=a$ implies $g(z)=b$, and {\it vice versa}. Sharing {\sc CM}
(counting multiplicities) means that $f-a$ ($1/f$ if $a=\infty$) and $g-b$ ($1/g$ if $b=\infty$) even have the same divisor, while {\sc IM}  (ignoring multiplicities)
means that nothing is assumed concerning multiplicities (the term `shared {\sc IM}' is established in the literature but just means `shared'). If $a=b$ we will just say that $f$ and $g$ share the {\it value} $a$.
The first result on functions sharing pairs of values is due {\sc Czubiak} and {\sc Gundersen}, who proved

{\sc Theorem A~ (\cite{CzubiakGundersen})}  {\it Meromorphic functions $f$ and $g$ sharing more than five pairs are M\"obius transformations of each other.}

Of course, functions that are M\"obius transformations of each other are {\it explicitly excluded} from consideration.

At the centre of the {\sc 4IM+1CM}-Problem is the question of whether or not the functions
\be\label{Gundersenex} \widehat f(z)=\frac{e^z+1}{(e^z-1)^2}\und \widehat g(z)=\frac{(e^z+1)^2}{8(e^z-1)}\ee
are essentially the only ones that share four pairs of values {\sc IM} and one pair {\sc CM}. They were constructed by {\sc Gundersen}~\cite{GGG3} as an example of functions that share four values, namely $0, \infty, 1, -1/8$, none of them {\sc CM}. Actually, the values $0$ and $1$ are assumed simply by $\widetilde f$ and doubly by $\widetilde g$, while for $-1/8$ and $\infty$ the reverse is true. Moreover, $\widehat f$ and $\widehat g$ share the pair $(-1/2,1/4)$ {\sc CM}. This was noticed by {\sc Reinders}~\cite{Reinders2}, who characterised the functions (\ref{Gundersenex})
in several ways.

{\sc Theorem B~ (\cite{Reinders3, Reinders2})} {\it Suppose meromorphic functions $f$ and $g$ share four values $a_\nu$ $(1\le\nu\le 4)$. Then
\be\label{Charakter}f=M \circ \widehat f\circ h\und g=M \circ \widehat g\circ h\ee
holds for some M\"obius transformation $M$ and some non-constant entire function $h$, provided one of the following conditions is fulfilled:

\begin{tabular}{rl}
{\rm (i)}&$f$ and $g$ share some extra pair $(a,b)$;\cr
{\rm (ii)}&for every $\nu\in\{1,2,3,4\}$, either $f-a_\nu$ or $g-a_\nu$ has only multiple zeros;\cr
{\rm (iii)}&for every $\nu\in\{1,2,3,4\}$, $(f-a_\nu)(g-a_\nu)$ has only triple zeros.
\end{tabular}}

The functions (\ref{Gundersenex}) also prove that the number five in Theorem A is best possible. In the context of functions sharing pairs of values, condition (\ref{Charakter})
has to be modified as follows to express the fact that $\widehat f$ and $\widehat g$ are {\it essentially unique:}
\be\label{fhatfghatg}f=M_1 \circ \widehat f\circ h\und g=M_2 \circ \widehat g\circ h.\ee
The next step after Theorem A was made by {\sc Gundersen}, {\sc Toghe}, and the author. They showed that at most one of five pairs of values can be shared {\sc CM}.

{\sc Theorem C~ (\cite{GST})} {\it Meromorphic functions $f$ and $g$ that are not M\"obius transformations of each other cannot share five pairs of values,
two of them {\sc CM}.}

Interestingly, the next Theorem was proved before Theorem C. On one hand it is stronger than Theorem C since only one {\sc CM}-pair is involved, and weaker on the other
by the additional hypothesis (\ref{Bedingung}) imposed on the proximity functions of the {\sc CM}-pair.

{\sc Theorem D~ (\cite{NSt5})} {\it Suppose that meromorphic functions $f$ and $g$ share five pairs of values. Under the additional hypothesis that one of these pairs,
$(a,b)$, say, is shared {\sc CM} and satisfies
\be\label{Bedingung}m(r,1/(f-a))+m(r,1/(g-b))=S(r),\ee
either $f$ and $g$ are M\"obius transformations of each other or else are given by} (\ref{fhatfghatg}).

The following best approximation to the corresponding {\sc 3IM+1CM}-Theorem for functions sharing four values is cited because it is closely related to Theorem~\ref{Coraltern}
in section~\ref{MAIN}.

{\sc Theorem E~ (\cite{GGG6})}  {\it Suppose that transcendental meromorphic functions $f$ and $g$ share three values {\sc IM} and one fourth value {\sc CM} with counting function
$N(r).$ Then either $N(r)\le \frac 45 T(r)+S(r)$ holds or else  $f$ and $g$ share all four values {\sc CM}.}

Condition (\ref{Bedingung}) was used in \cite{NSt5} to derive a quadratic relation between $f$ and $g$, that is, $f$ and $g$ parameterise some quadratic algebraic curve
(necessarily of genus zero).  In this paper it will be shown that also in the contrary case $f$ and $g$
parameterise some algebraic curve of genus zero of low degree (at most $9$) and depending on few (at most $7$) parameters. In this way the {\sc 3IM+1CM}-problem becomes a computer algebra problem.

Our main tool will be the theory of meromorphic functions and its interaction with uniformisation and parameterisation of algebraic curves. We will frequently make use of the {\sc Nevanlinna} functions $m(r,f), N(r,f), \overline N(r,f), T(r,f)$ and the First and Second Main Theorem of {\sc Nevanlinna} Theory. By $S(r,f)$ we denote any function
satisfying $S(r,f)=o(T(r,f))$ as $r\to\infty$, possibly outside some exceptional set of finite measure. For details the reader is referred to {\sc Hayman}'s book~\cite{Hayman}.
In the next section the most important properties of algebraic curves and some examples that are related to pair-sharing meromorphic functions will be discussed.

\section{Algebraic Curves}\label{ALGCURVES}

Besides the functions $e^z$ and $e^{-z}$, essentially three pairs of meromorphic functions $f$ and $g$ are known to share either four values or five pairs (Example~\ref{EX1} below due to {\sc Gundersen}~\cite{GGG3}, Example~\ref{EX4} due to {\sc Reinders}~\cite{Reinders1}, and Example~\ref{EX2} due to the author~\cite{NSt3-4})). In any case, $f$ and $g$ are algebraically dependent, that is,
they satisfy some non-trivial polynomial equation $K(f(z),g(z))=0$. The polynomial $K(x,y)$ of two complex variables may assumed irreducible.
Then the set
$$\mathcal K=\{(x,y)\in\C\times\C:K(x,y)=0\}$$
is called an {\it algebraic curve}. If it is necessary or desirable to add points like $(a,\infty)$ or $(\infty,b)$ or even $(\infty,\infty)$, one has to consider in addition
equations $K(a+x,1/y)y^m=0$ or $K(1/x,b+y)x^n=0$ or $K(1/x,1/y)x^ny^m=0$ at $(x,y)=(0,0)$; $m$ and $n$ denote the degree of $K$ w.r.t.\ $y$ and $x$, respectively.
The functions $f$ and $g$ are said to {\it parameterise} $\mathcal K$.

The following facts are taken from \cite{BeardonNg}. If $\mathcal K$ is parameterised by meromorphic functions in the plane, the curve has {\it genus} zero or one.

\medskip In case of genus {\it one}, any parameterisation has the form
\be\label{genus01}f=Q\circ h\und g=\widetilde Q\circ h,\ee
where $Q$ and $\widetilde Q$ are elliptic functions and $h$ is non-constant entire.

\medskip Much more important for us is the case of genus {\it zero}. By Theorem 1 and 4 in \cite{BeardonNg},
any parameterisation again has the form (\ref{genus01}), where now $Q$ and $\widetilde Q$ are rational functions and $h$ is any non-constant meromorphic on the plane. Moreover,
$$\mathcal K=\{(f(z),g(z)):z\in \C\}\cup {\mathcal E}$$
holds; the {\it exceptional set} ${\mathcal E}$ is finite and consists of {\it asymptotic values} $(a,b)$ of $(f,g)$, that is,
$(f(z),g(z))$ tends to $(a,b)$ as $z\to\infty$ on some plane curve. The rational functions $Q$ and $\tilde Q$ may be chosen such that the map
$\Phi=(Q,\widetilde Q)$ is injective on $\CD\setminus E$, where $E$ is a finite set (not to be confused with $\mathcal E$). With this choice, (\ref{genus01}) holds for every parametrisation
with suitably chosen meromorphic function $h$.

\begin{example}\label{EX0} \rm The algebraic curve $\mathcal C$ defined by $x^2+y^2=1$ $(x,y\in\C)$ has genus zero and the well-known rational and entire parameterisation
$$Q(t)=\ds\frac{1-t^2}{1+t^2},~\widetilde Q(t)=\ds\frac{2t}{1+t^2}\quad(t\ne\pm i)~\rm and$$
$$\ts f(z)=\cos z=Q(\tan \frac z2), ~g(z)=\sin z=\widetilde Q(\tan \frac z2),$$
respectively. On the other hand, $\overline{\mathcal C}=\mathcal C\cup\{\infty,\infty\}$ is parameterised over $\C$ by {\sc Jacobi}'s elliptic functions {\it sinus} and {\it cosinus amplitudinis}. This, however, does not mean that this curve has genus one but just says that the functions ${\sf cn}$ and $\sf sn$ may be written as $\ds\frac{1-h^2}{1+h^2}$ and $\ds\frac{2h}{1+h^2}$, respectively, where $h$ is some elliptic function (with different lattice) with simple zeros and poles at the zeros of $\sf sn$.
\end{example}

\begin{rem}\label{REM1} \rm In case of genus one it is obvious that it suffices to consider elliptic functions $f=Q$, $g=\widetilde Q$ in place of $f=Q\circ h$, $g=\widetilde Q\circ h$.
In case of genus zero we will often take the opportunity to switch from (\ref{genus01}) to functions
\be\label{genuszeroexp}f=Q\circ M\circ\exp\und g=\widetilde Q\circ M\circ\exp\ee
with some M\"obius transformation $M$. This is possible by the following reason: if $f$ and $g$ are given by (\ref{genus01}) and share five pairs of values, then $Q$ and $\widetilde Q$ also share these pairs on $\CD\setminus\{$Picard values of $h\}$.
In any case, if $h$ has either two Picard values ($a=M(\infty)$ and $b=M(0)$) or just one ($a=M(\infty)$) or none ($M={\rm id}$), also $Q\circ M\circ\exp\circ\phi$ and
$\widetilde Q\circ M\circ\exp\circ\phi$ for any non-constant entire function $\phi$ share the very same pairs {\sc IM/CM} on $\C$. The only thing that matters is the determination of $Q$ and $\widetilde Q$.
\end{rem}

\begin{example}\label{EX1} \rm The algebraic curve defined by the polynomial
\be\label{AlgEq} 4x^2+2cxy+y^2-8x\ee
is strongly related to the {\sc 4IM+1CM}-Problem. Like any quadratic curve it has genus zero. For $y=tx$ we obtain $x((4+2ct+t^2)x-8)=0$, hence the curve has the rational parameterisation
\be\label{para}x=Q(t)=\frac{8}{4+2ct+t^2},~y=\widetilde Q(t)=\frac{8t}{4+2ct+t^2}.\ee
Obviously, $Q$ and $\widetilde Q$ have poles at $t=-c\pm\sqrt{c^2-4}$; they are simple if $c\ne \pm 2$,
which together with $c\ne 0$ will henceforth be assumed. Hence the functions $f$ and $g$ given by (\ref{genus01})
share the value $\infty$ {\sc CM} for every non-constant meromorphic function $h$.
Suppose $f$ and $g$ also share the finite pairs $(a_\nu,b_\nu)$ ($1\le \nu\le 4$) {\sc IM} and {\it not} {\sc CM}. Then $a_\nu$ resp. $b_\nu$ is a critical value of $Q$ resp. $\widetilde Q$
and the following holds:

\begin{center}\begin{tabular}{c|c}
critical values&non-critical values\cr\hline
&\cr
$\begin{aligned}
Q(\infty)&=0\cr
Q(-c)&=8/(4-c^2)\cr
\widetilde Q(2)&=4/(2+c)\cr
\widetilde Q(-2)&=-4/(2-c)\end{aligned}$&
$\begin{aligned}
\widetilde Q(\infty)&=0=\widetilde Q(0)\cr
\widetilde Q(-c)&=-8c/(4-c^2)=\widetilde Q(-4/c)\cr
Q(2)&=2/(2+c)=Q(-2-2c)\cr
Q(-2)&=2/(2-c)=Q(2-2c)
\end{aligned}$
\end{tabular}\end{center}

If $c\ne 0,\pm 2$, $Q$ and $\widetilde Q$ share the value $\infty$ {\sc CM} and the pairs
\be\label{pairs}(0,0),~\Big(\frac 8{4-c^2},\frac{-8c}{4-c^2}\Big),~\Big(\frac 2{2+c}, \frac 4{2+c}\Big), {\rm ~and~} \Big(\frac 2{2-c},\frac{-4}{2-c}\Big)\ee
{\sc IM} on $\CD\setminus\{0,-4/c,-2-2c,2-2c\}$. Thus there is no non-constant meromorphic function $h$ in the plane such that the functions (\ref{genus01})
share the value $\infty$ and the pairs (\ref{pairs}) -- {\it except when} $c=\pm 1$. Then the points $-4/c$ and $\mp 2-2c$ and also $0$ and $\pm 2-2c$ coincide and the sphere is twice punctured at $-4,0$ and $4,0$, respectively. For any non-constant meromorphic function $h$ with Picard values  $-4$ and $0$  resp.\ $4$ and $0$, the functions (\ref{genus01}) share the value $\infty$ {\sc CM} and the pairs $(0,0), (8/3,-8/3), (2/3,4/3), (2,-4)$ resp. $(0,0), (8/3,8/3), (2,4), (2/3,-4/3)$. We note that then also (\ref{fhatfghatg}) holds.\end{example}

\begin{example}\label{EX4} \rm (\cite{Reinders1}) The algebraic curve defined by
$$(y-x)^4-16xy(x^2-1)(y^2-1)=0$$
has genus one. Every parameterisation $(f,g)$ shares the values $0,1,-1$ and $\infty$, with alternating multiplicities $(1\!:\!3)$ and $(3\!:\!1)$. For example, this means that $f$ has simple and triple zeros where $g$ has triple and simple zeros, respectively. In the most simple case, $f$ and $g$ are elliptic functions of elliptic order four.
\end{example}

\begin{example}\label{EX2} \rm (\cite{NSt3-4,NTNFADE}) The algebraic curve $\mathcal H$ defined by
\be\label{EXK}H(u,y)=y^3-3((\bar a-1)u^2-2u)y^2-3(2u^2-(a-1)u)y-u^3=0\ee
with $a=-\frac12+\frac i2\sqrt 3$ has genus zero and the rational parameterisation
$$u=\widetilde R(t)=3(a-1)\frac{t(1+t)^2}{(1+3t)^2},~ y=R(t)=3(a-1)\frac{t^2(1+t)}{1+3t}.$$
It is also parameterised by elliptic functions $U$ and $f$, where $U$ is a modified $\wp$-function satisfying $U'^2=U(U+1)(U-a)$.
The functions  $U$ and $f$ share the values $0$ and $\infty$ and the pairs $(a,1)$ and $(-1,-a)$, each with alternating multiplicities
$(2\!:\!4)$ and $(2\!:\!1)$. On its primitive lattice $\omega\Z+\omega'\Z$, $U$ has elliptic order two, while on their common lattice $(2\omega-\omega')\Z+(\omega+\omega')\Z$ both functions have elliptic order six. We note that at $(u,y)=(0,0)$, say, the equation $H(u,y)=0$ has the solutions $y\sim u^2$ (this meaning $y(u)=C_2u^2+C_3u^3+\cdots$, $C_2\ne 0$) and $u\sim y^2$, whose graphs parameterise $\mathcal H$ in a neighbourhood of $(0,0)$.
Actually there are three functions $f_1=f$, $f_2(z)=f(z+\omega)$, and $f_3(z)=f(z+\omega')$ of this kind. They share the values $0,\infty,1,-a$, each with alternating multiplicities $(4\!:\!1\!:\!1)$, $(1\!:\!4\!:\!1)$, and $(1\!:\!1\!:\!4)$. Up to transformations $f_\nu\mapsto M\circ f_\nu\circ h$ ($M$ any M\"obius transformation and $h$ any non-constant entire function) this triple is unambiguously determined by the requirement that three mutually distinct meromorphic functions share four values {\sc IM}.
Any two of these functions parameterise some algebraic curve $\mathcal K: K(x,y)=0$ of genus one. The polynomial $K$ is an appropriate factor of the resultant of the polynomials $K(u,y)$ and $K(u,x)$ with respect to the variable $u$; $\mathcal K$ is parameterised by elliptic functions $f_j=R\circ h_j$ and $f_k=R\circ h_k$ ($1\le j<k\le 3$) of elliptic order six; $h_j$ and $h_k$ are elliptic functions of elliptic order two w.r.t.\ the lattice $(2\omega-\omega')\Z+(\omega+\omega')\Z$ such that $\widetilde R\circ h_j=\widetilde R\circ h_k=U$ holds.
Note that $K(0,y)=y^6$, $K(1,y)=(y-1)^6,$ and $K(-a,y)=(y+a)^6$. At $(x,y)=(0,0)$, say, the equation $K(x,y)=0$ has solutions $y=-x+O(x^2)$, $y=cx^4+O(x^5)$, and $x=cy^4+O(y^5)$
with $c=i\sqrt 3/{243}$.
\end{example}

\section{Results}\label{MAIN}
From now on it is assumed that the following hypothesis holds:

\begin{itemize}\item[{\sc (H)}] $f$ and $g$ are transcendental meromorphic functions that are not M\"obius transformations of each other and share
four finite pairs $(a_\nu,b_\nu)$ {\sc IM}  and the value $\infty$ {\sc CM}, with counting functions $\overline N(r;a_\nu,b_\nu)$ and $\overline N(r)$, respectively.
\end{itemize}

We note that the definitions of sharing {\sc IM} and {\sc CM} may be relaxed insofar as
\begin{itemize}
\item[$\bullet$]zeros of $f-a_\nu$  that are not zeros of $g-b_\nu$, and {\it vice versa},
\item[$\bullet$] poles of $f$ that are not poles of $g$, and {\it vice versa}, and
\item[$\bullet$] poles of both functions with different multiplicities
\end{itemize}

are admitted, provided these points form sequences with counting function $S(r)$. These `generalisations' are so obvious that we don't need to say another word about it.

\begin{theo}\label{AlgDependent} Under hypothesis {\sc (H)}, $f$ and $g$ parameterise some algebraic curve $K(x,y)=0$ of genus zero.\end{theo}

The rather long and technical proof will be given in sections~\ref{sectA} and \ref{sect2B}.
The polynomial $K$ is found in semi-explicit form as divisor of a rather complicated polynomial of degree $9$ w.r.t.\ $x$ and $y$, and of degree $13$ w.r.t.\ $(x,y)$.
It turns out that this polynomial has to fulfill so many extra conditions that it is hard to believe that the {\sc 4IM+1CM}-Problem has solutions other that (\ref{fhatfghatg}).
In any case the following addendum to Theorem~\ref{AlgDependent} is valid and will certainly prove useful in future considerations.

\smallskip{\sc Addendum to Theorem~\ref{AlgDependent}.} {\it For $1\le\nu\le 4$,
\begin{itemize}
\item[I.] there exist integers $p_\nu>1$ and $q_\nu>1$ $(1\le\nu\le 4)$ with the following property: up to a sequence with
counting function $S(r)$, every solution of the equation $(f(z),g(z))=(a_\nu,b_\nu)$ has multiplicity $(p_\nu\!:\!1)$ or $(1\!:q_\nu)$ or $(1\!:\!1)$;
\item[II.] the polynomials $K(x,b_\nu)$ and $K(a_\nu,y)$ vanish at $x=a_\nu$ and $y=b_\nu$, respectively, and nowhere else in the plane.
\end{itemize}}

The proof will be given in Section~\ref{ADD}.

Our next result is similar to a combination of Theorem B~(iii) and Theorem E. It only formally contains Theorem D since this theorem will be part of the proof.

\begin{theo}\label{Coraltern} Suppose that in addition to hypothesis {\sc (H)} the zeros of $(f-a_\nu)(g-b_\nu)$ $(1\le\nu\le 4)$ have multiplicities at least three. Then either
\be\label{57Bed}\overline N(r)\le \frac 57 T(r)+S(r)\ee
holds or else the conclusion of Theorem D remains valid.\end{theo}

The multiplicities $(p_\nu\!:\!1)$ and $(1\!:\!q_\nu)$ may alternate, that is, the equation
$$(f(z),g(z))=(a_\nu,b_\nu)$$
may have $(p_\nu\!:\!1)$-fold and also $(1\!:\!q_\nu)$-fold solutions, not to mention solutions with multiplicity $(1\!:\!1)$. If, however, for each pair $(a_\nu,b_\nu)$ the multiplicity is always $(p_\nu\!:\!1)$ or always $(1\!:\!q_\nu)$, then
the desired {\sc 4IM+1CM}-Theorem can be deduced from Theorem D. This is the contents of the following theorem, which may be viewed as the analogue to Theorem B~(ii) for functions sharing pairs of values rather than values.

\begin{theo}\label{multiplicities1} Suppose {\sc (H)} holds and that for each $\nu$ either $f-a_\nu$ or else $g-b_\nu$ has only multiple zeros. Then the conclusion of Theorem D is valid.\end{theo}
Both theorems will be proved in section~\ref{Thm2and3}.

\begin{rem}\label{REM3} \rm We note that alternating multiplicities and also common simple zeros of $f-a_\nu$ and $g-b_\nu$ cannot be excluded {\it a priori} if only an algebraic curve serves as definition of the functions in question. This is shown by Examples~\ref{EX4} and \ref{EX2} in section~\ref{ALGCURVES}. We also note that {\it local} investigations of algebraic functions cannot utilise any information about the genus of the curve in question.\end{rem}

\section{Proof of Theorem~\ref{AlgDependent}: Preliminary Results}\label{sectA}

From various sources \cite{CzubiakGundersen, GGG6,GST,NSt5} one can deduce that up to the remainder term $S(r)$ of {\sc Nevanlinna} theory, meromorphic functions $f$ and $g$ that share
four pairs {\sc IM} and the value $\infty$ {\sc CM} have the same characteristic $T(r)$, proximity function of infinity $m(r)$, and counting function of poles $N(r)=\overline N(r)+S(r)$; the latter relation says that up to a sequence with counting function $S(r)$, the poles of $f$ and $g$ are simple.
Choose $(c_1,c_2,\ldots,c_5)\in\C^5$ non-trivially such that
$$P(x,y)=c_1x^2+c_2xy+c_3x+c_4y+c_5$$
vanishes at $(x,y)=(a_\nu,b_\nu)$, $1\le\nu\le 4$. Then
\be\label{PdefF}F(z)=P(f(z),g(z))\ee
does not vanish identically, since otherwise $g=-\ds\frac{c_1f^2+c_3f+c_5}{c_2f+c_4}$
would be a M\"obius transformation of $f$ (which is excluded) if $c_1=0$ or else would satisfy $T(r,g)=2T(r,f)+O(1)$ in contrast to $T(r,f)\sim T(r,g)$.
Then  the inequalities
\be\label{TrF}
\begin{aligned}
m(r,F)\le&~2m(r,f)+m(r,g)+S(r)=3m(r)+S(r)\cr
N(r,F)\le&~2\overline N(r)+S(r)\cr
T(r,F)\le&~2T(r)+m(r)+S(r)
\end{aligned}\ee
hold, and from
\be\label{sumN}\sum_{\nu=1}^4\overline N(r;a_\nu,b_\nu)\le \overline N(r,1/F)=T(r,F)-m(r,1/F)-N_1(r,1/F)+O(1)\ee
and the Second Main Theorem it follows that
$$\begin{aligned}
3T(r)\le &~\overline N(r)+{\ts\sum_{\nu=1}^4} \overline N(r;a_\nu,b_\nu)+S(r)\cr
\le &~\overline N(r)+T(r,F)-m(r,1/F)-N_1(r,1/F)+S(r)\cr
\le &~3T(r)-m(r,1/F)-N_1(r,1/F)+S(r)\cr
\le &~ 3T(r)+S(r).
\end{aligned}$$
Thus not only in this chain of inequalities the equality sign must hold everywhere but also in (\ref{TrF}) and (\ref{sumN}). In particular, this means that the poles and zeros of $F$ essentially arise from the poles and the finite shared pairs of values of $(f,g)$ and up to sequences with counting function $S(r)$ the zeros and poles of $F$ have multiplicity one and two, respectively. We have thus reached our first milestone.

\begin{lem}\label{C1} Under hypothesis {\sc (H)} the following holds:
\be\label{milestone}\begin{array}{rrcl}
{\rm (i)}&m(r,1/F)+N_1(r,1/F)&=&S(r)\cr
{\rm (ii)}&m(r,F)&=&3m(r)+S(r)\cr
{\rm (iii)}&N(r,F)&=&2\overline N(r)+S(r)\cr
&&=&2\overline N(r,F)+S(r)\cr
{\rm (iv)}&\sum_{\nu=1}^4 \overline N(r;a_\nu,b_\nu)&=&2T(r)+m(r)+S(r)\cr
&&=&\overline N(r,1/F)+S(r).
\end{array}\ee
\end{lem}
The same holds if $F$ is replaced by $\widetilde F(z)=\widetilde P(f(z),g(z))$, where
$$\widetilde P(x,y)=\widetilde c_1y^2+\widetilde c_2xy+\widetilde c_3x+\widetilde c_4y+\widetilde c_5$$
non-trivially satisfies $\widetilde P(a_\nu,b_\nu)=0$ ($1\le\nu\le 4)$.

For $c_1=0$ the better estimate $T(r,F)\le 2T(r)+S(r)$ is obtained, which implies $m(r)=S(r)$ and immediately gives (\ref{fhatfghatg}) by Theorem D.
The same is true if $\widetilde c_1=0$. For this reason we henceforth assume $c_1=\widetilde c_1=1$, hence
\be\label{PtildeP}\begin{aligned}
P(x,y)=&~ x^2+c_2xy+c_3x+c_4y+c_5 ~{\rm and}\cr
\widetilde P(x,y)=&~y^2+\widetilde c_2xy+\widetilde c_3x+\widetilde c_4y+\widetilde c_5.
\end{aligned}\ee
\begin{rem}\label{REM4} \rm
Various Polynomials of smallest possible degree vanishing at the shared pairs of values were successfully used
in \cite{CzubiakGundersen, GGG6,GST,NSt5}. The quadratic polynomials in $(x,y)$ that vanish at the points $(a_\nu,b_\nu)$ form a linear space.
The polynomials (\ref{PtildeP}) form a basis except when some nontrivial polynomial
$$P_0(x,y)=\widehat c_2xy+\widehat c_3x+\widehat c_4 y+\widehat c_5$$
also vanishes at each $(a_\nu,b_\nu)$. Then $b_\nu=M(a_\nu)$ holds for some M\"obius transformation
$M$, $f$ and $M^{-1}\circ g$ share four values $a_\nu$ and the pair $(\infty,M^{-1}(\infty))$ and Theorem B applies.
This is the case for the pairs $(0,0)$, $(8/3,-8/3)$, $(2,-4)$, $(2/3,4/3)$; $P_0(x,y)=3xy+4x-4y$ and $\widetilde P_0(x,y)=12x^2+3y^2-32x+8y$  form a basis,
and the polynomials (\ref{PtildeP}) are not available, that is, $c_1=\widetilde c_1=0$ and $P$ and $\widetilde P$ are constant multiples of $P_0(x,y)=3xy+4x-4y$.
\end{rem}

The {\it Proof} of Theorem~\ref{AlgDependent} requires a sequence of rather technical results which are based on Lemma~\ref{C1} and the hypothesis
\be\label{mneS}m(r)\ne S(r).\ee
From now on this will be assumed
in addition to hypothesis {\sc (H)}. Since the author believes that (\ref{mneS}) is fictional, various auxiliary results will be called {\sc Claim} rather than {\sc Lemma}. First of all we will show that  $|f|$ and $|g|$ are `large' on disjoint sets.
To be more precise, set
$$\begin{aligned}
E_r=&~\{\theta\in (-\pi,\pi]:|f(re^{i\theta})|\ge|g(re^{i\theta})|\}\quad{\rm and}\cr
E_r^*=&~\{\theta\in (-\pi,\pi]:|f(re^{i\theta})|\ge K(1+|g(re^{i\theta})|^2)\}\subset E_r;
\end{aligned}$$
$K>1$ will be fixed later. The sets $\widetilde E_r$ and $\widetilde E_r^*$ are defined in the same way with $f$ and $g$ interchanged.

\begin{claim}\label{C2} $\ds m(r)=\ds\frac 1{2\pi}\int_{E_r^*}\log^+|f|\,d\theta+S(r)=\frac 1{2\pi}\int_{\widetilde E_r^*}\log^+|g|\,d\theta+S(r)$.\end{claim}

\proof The inequalities (for $F(z)=P(f(z),g(z))$)
$$\log^+|F|\le 2\log^+|f|+O(1)\quad{\rm on~} E_r$$
(following from $|g|\le |f|$) and
$$\log^+|F|\le \log^+|f|+\log^+|g|+O(1)\quad{\rm on}~\widetilde E_r$$
(use $|f^2|\le |f||g|$) imply
$$\begin{aligned}m(r,F)&\le\frac 2{2\pi}\int_{E_r}\log^+|f|\,d\theta+\frac 1{2\pi}\int_{\widetilde E_r}\log^+|f|\,d\theta+\frac 1{2\pi}\int_{\widetilde E_r}\log^+|g|\,d\theta+O(1)\cr
&\le 3m(r)-\frac 1{2\pi}\int_{\widetilde E_r}\log^+|f|\,d\theta-\frac 1{2\pi}\int_{E_r}\log^+|g|\,d\theta+O(1).\end{aligned}$$
Since, however, $m(r,F)=3m(r)+S(r)$ holds by (\ref{milestone})(ii), this implies
$$\begin{aligned}
m(r,f)=&~\frac 1{2\pi}\int_{E_r}\log^+|f|\,d\theta+S(r),\cr
m(r,g)=&~\frac 1{2\pi}\int_{\widetilde E_r}\log^+|g|\,d\theta+S(r),\quad{\rm and}\end{aligned}$$
$$\frac 1{2\pi}\int_{\widetilde E_r}\log^+|f|\,d\theta+\frac 1{2\pi}\int_{E_r}\log^+|g|\,d\theta=S(r).$$
The assertion then follows from
$$\frac1{2\pi}\int_{E_r\setminus E_r^*}\log^+|f|\,d\theta\le \log K+\frac1{2\pi}\int_{E_r}\log^+(1+|g|^2)\,d\theta=S(r).\eqno{\ende}$$

\begin{claim}\label{phitildephi} The functions $\ds \phi=\frac{f'F^2/\widetilde F}{\prod_{\nu=1}^4(f-a_\nu)}$ and $\ds\widetilde \phi=\frac{g'\widetilde F^2/F}{\prod_{\nu=1}^4(g-b_\nu)}$
have {\sc Nevanlinna} characteristic $S(r)$.
\end{claim}

\proof Up to a sequence of points with counting function $S(r)$, poles of $\phi$ may only occur at multiple zeros  of $\widetilde F$  and at simple  poles of $\widetilde F$. Since, however, again
up to sequences with counting function $S(r)$, the zeros of $\widetilde F$ are simple  by (\ref{milestone})(i) (for $\widetilde F$ in place of $F$) and the poles of $F$ and $\widetilde F$ have order two, again by (\ref{milestone})(iii) (for $\widetilde F$ and $F$)
we have $N(r,\phi)\le N(r,F/\widetilde F)+N(r,\widetilde F/F)=S(r)$. To prove $m(r,\phi)=S(r)$ we set
$$L=\frac{f'}{\prod_{\nu=1}^4(f-a_\nu)}$$
and note that $m(r,1/\widetilde F)=S(r)$ holds by (\ref{milestone})(i) (for $\widetilde F$ in place of $F$).
To estimate $m(r,\phi)$ we consider the contributions of the mutually disjoint sets
$E_r^*$,  $\widetilde E_r^*$, and  $\widehat E_r=(-\pi,\pi]\setminus (E_r^*\cup \widetilde E_r^*)$, noting that
$$\frac1{2\pi}\int_{E_r^*\cup \widehat E_r}\log^+|g|\,d\theta+\frac1{2\pi}\int_{\widetilde E_r^*\cup \widehat E_r}\log^+|f|\,d\theta=S(r).$$

On $\widetilde E_r^*$ (where $|g|\ge K(1+|f|^2)$ is `large') we have $|F|^2=O((1+|f|^2)^2|g|^2)$ and $|\widetilde F|>\frac 12|g|^2$
if $K$ is chosen sufficiently large, hence $|F|^2/|\widetilde F|=O((1+|f|^2)^2)$ and
the contribution of $\widetilde E_r^*$ to $m(r,\phi)$ is at most
$$m(r,L)+\frac1{2\pi}\int_{\widetilde E_r^*}\log\;(1+|f|^2)^2\,d\theta+O(1)=S(r).$$

On $\widehat E_r$ we have $|F|=O(|f|^2+|g|^2+1)$, hence the contribution of $\widehat E_r$ is at most
$$\frac 1{2\pi}\int_{\widehat E_r}\log(1+|f|^2+|g|^2)\,d\theta+m(r,L)+m(r,1/\widetilde F)+O(1)=S(r).$$
On $E_r^*$ write
\be\label{gdurchf}\phi=\frac{f'(f^2+c_2fg+c_3f+c_4g+c_5)^2}{(f-a)\prod_{\nu=1}^4(f-a_\nu)}\times\frac{(f-a)}{\widetilde F},\ee
where $a\ne a_\nu$ ($1\le\nu\le 4$) is any complex number.
Since $|f|\ge|g|$ holds on $E_r^*$ we may write
$$f^2+c_2fg+c_3f+c_4g+c_5=c^*_2f^2+c^*_3f+c_5,$$
where $c_2^*=1+c_2(g/f)$ and $c^*_3=c_3+c_4(g/f)$ are considered as coefficients satisfying $|c^*_2|\le 1+|c_2|$ and $|c^*_3|\le|c_3|+|c_4|$.
Thus the first factor in (\ref{gdurchf}) is
$$\ds O\Big(\frac{|f'|}{|f-a|}+\sum_{\nu=1}^4\frac{|f'|}{|f-a_\nu|}\Big)$$
(partial fraction decomposition) and has proximity function $S(r)$.
The First Main Theorem yields
$$m\Big(r,\frac{f-a}{\widetilde F}\Big)=m\Big(r,\frac{\widetilde F}{f-a}\Big)+N\Big(r,\frac{\widetilde F}{f-a}\Big)-N\Big(r,\frac{f-a}{\widetilde F}\Big)+O(1).$$
To estimate the first term to the right we note that on $E_r^*$ we have
$$\ds\Big|\frac{\widetilde F}{f-a}\Big|=O(1+|g|),$$
while $|\widetilde F|=O(|g|^2)$ and $|\widetilde F|=O(1+|f|^2+|g|^2)$ holds on $\widetilde E_r^*$ and $\widehat E_r$, respectively. This yields the estimate
 $$\begin{aligned}
m\Big(r,\frac{\widetilde F}{f-a}\Big)\le&~ m\Big(r,\frac 1{f-a}\Big)+\frac1{2\pi}\int_{E^*_r}\log\;(1+|g|)\,d\theta+\frac 2{2\pi}\int_{\widetilde E^*_r}\log^+|g|\,d\theta\cr
&~+\frac1{2\pi}\int_{\widehat E_r}\log\;(1+|f|^2+|g|^2)\,d\theta +O(1)\cr
=&~m\Big(r,\frac 1{f-a}\Big)+2m(r)+S(r).\end{aligned}$$
Combining this with
$$\ds N\Big(r,\frac{f-a}{\widetilde F}\Big)\ge N(r,1/\widetilde F)+S(r)=2T(r)+m(r)+S(r).$$
and
$$\ds N\Big(r,\frac{\widetilde F}{f-a}\Big)\le N\Big(r,\frac 1{f-a}\Big)+\overline N(r)+S(r)$$
yields
$$\begin{aligned}
m\Big(r,\frac{f-a}{\widetilde F}\Big)\le &~ 2m(r)+m\Big(r,\frac 1{f-a}\Big)+N\Big(r,\frac 1{f-a}\Big)+\overline N(r)\cr
&~-(2T(r)+m(r))+S(r)=S(r)\end{aligned}$$
and eventually $T(r,\phi)=S(r)$.\hfill\Ende

\begin{claim}\label{C4} $m(r,F/\widetilde F)=m(r)+S(r)$.\end{claim}

\proof Set $\Psi=F/\widetilde F$. To estimate $m(r)=m(r,f)+S(r)$ write
$$\ds\frac{f}\Psi=\frac{f\widetilde F}{F}=\frac{fg^2+\widetilde c_2f^2g+\widetilde c_3f^2+\widetilde c_4fg+\widetilde c_5f}{f^2+c_2fg+c_3f+c_4g+c_5}.$$
On $E_r^*$, where $|f|\ge K(1+|g|^2)$ is `large' and thus
$$\begin{aligned}
|f^2+c_2fg+c_3f+c_4g+c_5|&\ge \frac 12|f|^2\und\cr
|fg^2+\widetilde c_2f^2g+\widetilde c_3f^2+\widetilde c_4fg+\widetilde c_5f|&=O(|f|^2(1+|g|)
\end{aligned}$$
holds we have
$|f|=O((1+|g|)|\Psi|)$, hence
$$\begin{aligned}
m(r)&=\frac1{2\pi}\int_{E_r^*}\log^+|f|\,d\theta+S(r)\cr
&\le \frac1{2\pi}\int_{E_r^*}\log^+|\Psi|\,d\theta+\frac1{2\pi}\int_{E_r^*}\log(1+|g|)\,d\theta+S(r)\cr
&\le m(r,\Psi)+S(r).
\end{aligned}$$
Besides $L$ we will also consider $\ds\widetilde L=\frac{g'}{\prod_{\nu=1}^4(g-b_\nu)}.$
Then $L$ and $\widetilde L$ have proximity function $S(r)$ and $\ds\frac{\phi}{\widetilde \phi}=\frac{L}{\widetilde L}\Psi^3$ even has characteristic $S(r)$, hence
$$\begin{aligned}
3m(r,\Psi)&\le m(r,1/L)+m(r,\widetilde L)+T(r,\phi/\widetilde \phi)+S(r)\cr
&=N(r,L)-N(r,1/L)+S(r)\cr
&\le\ts\sum_{\nu=1}^4\overline N(r;a_\nu,b_\nu)-2\overline N(r)+S(r)\cr
&=2T(r)+m(r)-2\overline N(r)+S(r)\cr
&=3m(r)+S(r)\end{aligned}$$
holds. Altogether this is nothing but
$m(r)=m(r,\Psi)+S(r)=m(r,F/\widetilde F)+S(r)$ and
$\ds m(r)=\frac1{2\pi}\int_{E_r^*}\log^+|\Psi|\,d\theta+S(r)=\frac1{2\pi}\int_{\widetilde E_r^*}\log^+|1/\Psi|\,d\theta+S(r).$
\hfill\Ende

From this (and the corresponding resut for $\widetilde F/F$) our next claim follows immediately.

\begin{claim}$\psi=\Psi'/\Psi=F'/F-\widetilde F'/\widetilde F$ has {\sc Nevanlinna} characteristic $S(r)$.\end{claim}

We note that $\psi$ does not vanish identically since otherwise $\Psi=F/\widetilde F$ would be a constant in contrast to $m(r,\Psi)=m(r)+S(r)\ne S(r).$

\section{Proof of Theorem~\ref{AlgDependent}: Construction of an Algebraic Curve}\label{sect2B}

In
\be\label{DEFHVorform}\psi=\big({P_x}/{P} -{\widetilde P_x}/{\widetilde P}\big)f'+\big({P_y}/{P} -{\widetilde P_y}/{\widetilde P}\big)g',\ee
wherein $P$, $\widetilde P$, $P_x$, $P_y$, $\widetilde P_x$, $\widetilde P_y$ have to be evaluated at $(f(z),g(z))$, the derivatives $f'$ and $g'$ may be replaced by
$$f'=\phi\frac{\widetilde P}{P^2}{\prod_{\nu=1}^4(f-a_\nu)}\und g'=\widetilde \phi\frac{P}{\widetilde P^2}{\prod_{\nu=1}^4(g-b_\nu)},$$
respectively, thus $(x,y)=(f(z),g(z))$ satisfies some algebraic equation
$$H(z,x,y)=0$$
over the field generated by the functions  $\psi,\phi$, and $\widetilde\phi$. Actually, $H$ is given by
\be\label{DEFH}\psi P^3\widetilde P^3-\phi(\widetilde P P_x-P\widetilde P_x)\widetilde P^3\prod_{\nu=1}^4(f-a_\nu)+\widetilde\phi(\widetilde P P_y-P\widetilde P_y)P^3\prod_{\nu=1}^4(g-b_\nu).\ee

Our next task is to prove

\begin{claim}\label{Hnontrivial} $H$ is non-trivial.\end{claim}

\proof Suppose $H$ is trivial. Since $\psi(z)\not\equiv 0$, the product $P(x,y)\widetilde P(x,y)$ vanishes at every point $(a_\mu,b_\nu)$. We may assume $P(a_\mu,b_\nu)=0$ for $\mu=1,2$, say, and two values
of $\nu$ depending on $\mu$ (not necessarily the same in both cases). Then
$$P(a_\mu,b_\nu)=a_\mu^2+c_3 a_\mu+c_5+(c_2a_\mu+c_4)b_\nu=0$$
holds for two different values of $\nu$, and this implies $a_\mu^2+c_3 a_\mu+c_5=c_2a_\mu+c_4=0$ for $\mu=1,2$,
hence $c_2=c_4=0$ since $a_1\ne a_2$. Of course, $P(x,y)=x^2+c_3x+c_5$ cannot vanish at $(a_\nu,b_\nu)$ for $1\le \nu\le 4$, this showing that $H$ is non-trivial.\hfill\Ende

To get rid of the functions $\psi,\phi,\widetilde\phi$ we now have to attack the most involved part of the proof:

\begin{claim}\label{C7} The ratios $\phi/\psi$ and $\widetilde\phi/\psi$ are constant.\end{claim}

\proof The polynomial $H$ has degree at most degree nine with respect to the single variables $x$ and $y$, thus
$$H(z,x,y)=\sum_{j=0}^9h_j(z,y)x^j$$
holds. The assertion is trivially true if $h_9(z,x)$ vanishes identically: for $a_1=b_1=0$, say, $h_9$ is given (thanks to {\texttt{maple}}) by
$$\begin{aligned}
h_9(z,y)=&~\widetilde c_3^3(\psi(z)-\widetilde c_3\phi(z))+\widetilde c_2y\big\{3\widetilde c_3^2\psi(z)-4\widetilde c_3^3\phi(z)+b_2b_3b_4\widetilde\phi(z))\cr
&~+(3\widetilde c_2\widetilde c_3\psi(z)-6\widetilde c_2\widetilde c_3^2\phi(z)-(b_2b_3+b_3b_4+b_4b_2)\widetilde\phi(z))y\cr
&~+(\widetilde c_2^2\psi(z)-4\widetilde c_2^2\widetilde c_3\phi(z)+(b_2+b_3+b_4)\widetilde\phi(z))y^2-(\widetilde c_2^3\phi(z)+\widetilde\phi(z))y^3\big\}.\end{aligned}$$
Assuming $h_9(z,y)\equiv 0$, either
\be\label{psitildec3phi}\psi=\widetilde c_3\phi\quad({\rm and~}\widetilde c_3\ne 0)\ee
holds or else $\widetilde c_3=0$. In this case, $\widetilde\phi(z)\not\equiv 0$ implies $\widetilde c_2=0$ and $\widetilde P(x,y)=y^2+\widetilde c_4y$, which, of course, is impossible
and proves (\ref{psitildec3phi}).

We will now show that exactly the same is true if $h_9(z,y)\not\equiv 0$. From $H(z,x,y)=0$ and $|x|>1$ it follows that
$$|h_9(z,y)||x|\le \sum_{j=0}^8 |h_j(z,y)|=O(\max\{|\psi(z)|,|\phi(z)|,|\widetilde\phi(z)|\}(1+|y|)^9),$$
hence (set $x=f(z)$, $y=g(z)$)
$$\begin{aligned}
m(r)=&~\frac1{2\pi}\int_{E_r^*}\log|f(re^{i\theta})|\,d\theta+S(r)\cr
\le&~ m(r,1/h_9(z,g(z)))+ m(r,\psi)+ m(r,\phi)+ m(r,\widetilde\phi)+O(1)\cr
&~+\frac 9{2\pi}\int_{E_r^*}\log (1+|g(re^{i\theta})|)\,d\theta+S(r)\cr
=&~m(r,1/h_9(z,g(z)))+S(r).\end{aligned}$$
To proceed we need the following

\begin{lem}For at least one index $j$, $h_9(z,b_j)$ vanishes identically.\end{lem}

\proof Suppose to the contrary that $h_9(z,b_j)\not\equiv 0$ for $1\le j\le 4$. Then also $h(z,b_j)\not\equiv 0$ holds for every prime factor of $h_9$ and
$$G(z,y)=h(z,y)\prod_{\nu=1}^4(y-b_\nu)$$
is a square-free polynomial over the algebraic closure of the field of meromorphic functions with {\sc Nevanlinna} characteristic $O(T(r,\phi)+T(r,\widetilde\phi)+T(r,\psi))$.
Then from Corollary 4 in {\sc Yamanoi}'s paper \cite{Yamanoi} it immediately follows that
$$(\deg_y G-1)T(r)\le \overline N(r)+\overline N(r,1/G(z,g(z)))+\epsilon T(r)~\|$$
holds for every $\epsilon>0$; the symbol $\|$ means {\it outside some exceptional set that depends on $\epsilon$ and has finite measure with respect to $d\log\log r$}.
Since, however, already
$$\overline N(r)+\sum_{\nu=1}^4 \overline N(r,1/(g-b_\nu))=3T(r)+S(r)~\|$$
holds, this yields
$$\overline N(r,1/h(z,g(z))\ge (\deg_y h-\epsilon)T(r)+S(r)~\|.$$
Taking into account all prime factors according to their multiplicities we eventually  obtain
$$\overline N(r,1/h_9(z,g(z))\ge (\deg_y h_9-4\epsilon)T(r)+S(r)~\|,$$
hence
$$m(r)\le m(r,1/h_9(z,g(z)))+S(r)\le 4\epsilon T(r)+S(r)~{\|}$$
Although this is weaker than the condition $m(r)=S(r)$, the proof in ~\cite{NSt5} goes through in the very same way, that is, Theorem D again
holds. This eventually proves $h_9(z,b_j)\equiv 0$ for some $j$.\hfill\Ende

To simplify notation we may assume $j=1$ and $a_1=b_1=0$, hence $\widetilde c_5=0$ and
$$h_9(z,0)=\widetilde c_3^3(\psi(z)-\widetilde c_3\phi(z))$$
vanishes identically. In case of $\widetilde c_3=0$ (and $\psi-\widetilde c_3\phi\not\equiv 0$) we obtain
$$\widetilde P(x,y)=y(y+\widetilde c_2x+\widetilde c_4),$$
hence $b_\nu=-\widetilde c_2 a_\nu-\widetilde c_4=M(a_\nu)$ holds for $2\le\nu\le 4$ and also $\nu=5$ with $a_5=b_5=\infty$.
In other words, the functions $f$ and $M^{-1}\circ g$ share four values $a_\nu$ ($2\le\nu\le 5$) and the pair $(0,M^{-1}(0))$.
By Theorem B we are done: the goal, the representation (\ref{fhatfghatg}) has been arrived head of time.
This eventually proves (\ref{psitildec3phi}). The corresponding assertion $\psi=c_3\widetilde \phi$ is proved in exactly the same manner.\hfill\Ende

From (\ref{DEFH}) and Claim~\ref{C7} it then follows that $H(z,x,y)=\psi(z)H_0(x,y),$
hence also $H_0(f(z),g(z))=0$ holds. Of course, $H_0$ does not have to be irreducible, nevertheless we have

\begin{claim}\label{C8} $(f,g)$ parameterises an algebraic curve $\mathcal K:K(x,y)=0$ of genus zero.\end{claim}

To prove {\sc Claim}~\ref{C8}, which marks the end of proof of Theorem~\ref{AlgDependent}, we assume to the contrary that $\mathcal K$ has genus one and
$Q$ and $\widetilde Q$ in (\ref{genus01}) are elliptic functions. Then also $Q$ and $\widetilde Q$ share the pairs $(a_\nu,b_\nu)$ {\sc IM} (and not {\sc CM})
and the value $\infty$ {\sc CM}. Since
$$m(r,Q)+m(r,\widetilde Q)=O(1),$$
Theorem D applies. It yields $Q=M\circ \widehat f\circ k$ and $\widetilde Q=\widetilde M\circ \widehat g\circ k$ with M\"obius transformations $M$
and $\widetilde M$ and some non-constant entire function $k$. It is, however, obvious that $\widehat f\circ k$ and $\widehat g\circ k$ are not elliptic functions.
Thus the genus of our algebraic curve is zero as was stated in Theorem~\ref{AlgDependent}, and (\ref{genus01}) holds with rational functions $Q$ and $\widetilde Q$ and some meromorphic function $h$.\hfill\Ende

\begin{rem}\label{REM6}\rm The proof of Theorem~\ref{AlgDependent} is much easier for {\sc Yosida} functions $f$ and $g$, that is, for functions with bounded spherical derivative $f^\sharp=\ds\frac{|f'|}{1+|f|^2}$. In this case the functions $\phi,\widetilde\phi,\psi$ are constants and there is no need for Claim~\ref{C7} and its sophisticated
and elaborate proof. The first idea therefore was to apply the  re-scaling technique due to {\sc Zalcman}~\cite{Zalcman} simultaneously (like in \cite{NStVier}) to obtain {\sc Yosida} functions
$\widetilde f(z)=\lim_{n\to\infty}f(z_n+\rho_n z)$ and $\widetilde g(z)=\lim_{n\to\infty}g(z_n+\rho_n z)$
for suitably chosen sequences $z_n\to\infty$ and $\rho_n\to 0$. By {\sc Hurwitz}' Theorem, these functions share the same values and pairs of values as do $f$ and $g$. However, it cannot
be excluded that $\widetilde f$ and $\widetilde g$ are M\"obius transformations of each other although $f$ and $g$ are not. The same problem arises if the re-scaling technique is simultaneously applied to functions $f$ and $g$ that share four values $a_\nu$: while it is true that $\widetilde f$ and $\widetilde g$ share the same values as $f$ and $g$, $\widetilde f=\widetilde g$ cannot be excluded although $f\ne g$.\end{rem}

\section{Proof of the Addendum to Theorem~\ref{AlgDependent}}\label{ADD}

I. For the sake of simplicity let us consider the case $\nu=1$ and $a_1=b_1=0$. Assuming $(f(z_0),g(z_0))=(0,0)$ with multiplicity $(1\!:\!q)$ and $q>1$ yields $F'(z_0)=c_3f'(z_0)$ and
$\widetilde F'(z_0)=\widetilde c_3f'(z_0)$, hence
$$ \phi(z_0)=\frac{({c_3^2}/{\widetilde c_3})f'(z_0)}{\prod_{\mu=2}^4(-a_\mu)}\und \widetilde\phi(z_0)=\frac{q({\widetilde c_3^2}/{c_3})f'(z_0)}{\prod_{\mu=2}^4(-b_\mu)}.$$
Since, however, $\widetilde\phi/\phi$ is constant by {\sc Claim}~\ref{C7}, $q$ is independent of $z_0$.

II. Since only properties of the polynomial $K(x,y)$ are affected we may assume that our algebraic curve $\mathcal K: K(x,y)=0$ is parameterised by the functions
\be\label{fgExpo} f(z)=Q(e^z)\und g(z)=\widetilde Q(e^z),\ee
which share the value $\infty$ {\sc CM} and the pairs $(a_\nu,b_\nu)$ {\sc IM} (and not {\sc CM}). Note that  the rational functions $Q$ and $\widetilde Q$ share $\infty$ {\sc CM} and the pairs
$(a_\nu,b_\nu)$ {\sc IM} (and not {\sc CM}) certainly on $\C\setminus\{0\}.$
Also, since $m(r,f)\sim m(r,g)\sim m(r)\ne o(r)$ and $f$ and $g$ are `large' on disjoint sets, $Q$ and $\widetilde Q$ have a pole of the same order
at $t=\infty$ and $t=0$, respectively, or {\it vice versa}. The exceptional set $\mathcal E$ in
$$\mathcal K=\{(f(z),g(z)):z\in\C\}\cup \mathcal E$$
consists of the asymptotic values
\be\label{asympval}\begin{aligned}
\lim_{\xi\to +\infty}(Q(e^\xi),\widetilde Q(e^\xi))=&~(Q(\infty),\widetilde Q(\infty))=(\infty,b_\kappa)\und\cr
\lim_{\xi\to -\infty}(Q(e^\xi),\widetilde Q(e^\xi))=&~(Q(0),\widetilde Q(0))=(a_\lambda,\infty)\end{aligned}\ee
for some $1\le \kappa,\lambda\le 4$. Thus $x=\widetilde Q(t_0)=a_\nu$ for some $t_0\in\C\setminus\{0\}$ implies $y=Q(t_0)=b_\nu$, that is, $K(a_\nu,y)=0$ implies $y=b_\nu$
(with multiplicity $1$ or $q_\nu$). \hfill\Ende

\section{Proof of Theorems~\ref{Coraltern} and \ref{multiplicities1}}\label{Thm2and3}

To prove Theorem~\ref{Coraltern} we need the estimate
\be\label{multiplepoints}\overline N_1(r;1/(f-a_\nu))\le m(r)+S(r)\ee
and, of course, also $\overline N_1(r;1/(g-b_\nu))\le m(r)+S(r).$
To simplify notation we may assume $\nu=1$ and $a_1=b_1=0$, choose
$$R(x,y)=Ax^2+By^2+Cxy+Dx$$
non-trivially such that $R(a_\nu,b_\nu)=0$ holds for $2\le\nu\le 4$ and set
$$h(z)=R(f(z),g(z)).$$
Then $h$ vanishes at the $(a_\nu,b_\nu)$-points of $(f,g)$, at least twice at the multiple zeros of $f$ (with counting function $\overline N_1(r,1/f)$).
If $h$ does not vanish identically the Second Main Theorem again yields
$$\begin{aligned}
2T(r)+m(r)\le&~{\ts\sum_{\nu=1}^4}\overline N(r;a_\nu,b_\nu)+S(r)\cr
\le &~ N(r,1/h)-\overline N_1(r,1/f)+S(r)\cr
\le &~ T(r,h)-\overline N_1(r,1/f)+S(r)\cr
\le &~2N(r)+4m(r)-\overline N_1(r,1/f)+S(r)\cr
= &~2T(r)+2m(r)-\overline N_1(r,1/f)+S(r),\end{aligned}$$
hence (\ref{multiplepoints}) for $a_\nu=0$.
If, however, $h(z)\equiv 0$, the functions $f$ and $g$ parameterise the algebraic curve $\mathcal R: Ax^2+By^2+Cxy+Dx=0$.
By Remark~\ref{REM4} we have $AB\ne 0$ and also $D\ne 0$ since the algebraic curve $Ax^2+By^2+Cxy=0$ is reducible. By the elementary change of variables
$(x,y)\mapsto (\gamma\alpha x,\gamma\beta y)$, the equation $R(x,y)=0$ may be transformed into
$$\gamma^2[4x^2+2cxy+y^2-8x]=0$$
with $A\alpha^2=4, B\beta^2=1,C\alpha\beta=2c,$ and $D\alpha/\gamma=-8.$
By Example~\ref{EX1} it is necessary that $c=\pm 1$ and $f$ and $g$ are given by (\ref{fhatfghatg}).
Thus if Theorem~D does not hold and $f-a_\nu$ and $g-b_\nu$ have no common simple zeros,
$$\overline N(r;a_\nu,b_\nu)\le \overline N_1(r,1/(f-a_\nu))+\overline N_1(r,1/(g-b_\nu))+S(r)\le 2m(r)+S(r)$$
holds for $1\le\nu\le 4$, and so
$$2T(r)+m(r)\le \sum_{\nu=1}^4 \overline N(r;a_\nu,b_\nu)\le 8m(r)+S(r)$$
by (\ref{milestone})(iv). This implies (\ref{57Bed}).\hfill\Ende

{\it Proof} of Theorem~\ref{multiplicities1}. For $\nu$ fixed, the hypothesis on the multiplicities implies
$$\overline N(r;a_\nu,b_\nu)\le \frac 12T(r)+O(1)$$
in any case, hence $\sum_{\nu=1}^4\overline N(r;a_\nu,b_\nu)\le 2T(r)+O(1)$. Then
(\ref{milestone})(iv) yields $m(r)=S(r)$ and Theorem D applies.\hfill\Ende

\section{Why the {\sc 4IM+1CM}-Conjecture is Probably True}

The polynomial in Theorem~\ref{AlgDependent} has the form
\be\label{KstA}K(x,y)=(x-a_\lambda)^sy^m+A(y-b_\kappa)^tx^n+\sum_{j,k}c_{jk}x^jy^k\ee
and satisfies
$K(a_\nu,y)=0\Leftrightarrow y=b_\nu$ and $K(x,b_\nu)=0\Leftrightarrow x=a_\nu$ $(1\le\nu\le 4),$
hence $K(a_\nu,b_\nu+y)$ and $K(a_\nu+x,b_\nu)$ are monomials. This requires
\be\label{Cond1}\frac{\partial^j K}{\partial y^j}(a_\nu,b_\nu)=\frac{\partial^\ell K}{\partial x^\ell}(a_\nu,b_\nu)=0\ee
for all but one $j$ and $\ell$, respectively; $m\le 9$ and $n\le 9$ are the degrees of $K$ w.r.t.\ $y$ and $x$, respectively, $K$ has degree at most $13$ w.r.t.\ $(x,y)$, and $1\le s,t\le 4$ holds.
On the other hand, since two of the pairs $(a_\nu,b_\nu)$ may be prescribed ($a_1=b_1=0$ and $a_2=b_2=1$, say) and $\phi/\psi$ and $\widetilde\phi/\psi$ may be expressed in terms of the coefficients $c_j,\widetilde c_k$, there are only {\it five} free parameters $a_3,b_3,a_4,b_4,A$
at hand to satisfy the $4(n+m-1)$ constraints (\ref{Cond1}), not to mention the fact that $K$ has to be irreducible of genus zero.
Of course, counting algebraic equations and variables does not disprove the existence of $K$ since for special pairs $(a_\nu,b_\nu)$ some of the equations may be algebraically dependent.
Nevertheless it might well be that the {\sc 4IM+1CM} problem will be solved by a virtuoso in computer algebra systems rather than a complex analyst.

\section{\sc Concluding Remark} It is not implausible to believe that all pairs $(f,g)$ of meromorphic functions sharing four values or five pairs are already known in essence and are presented in Examples \ref{EX1}--\ref{EX2}. A proof of this very general assertion, however, seems to be beyond the present knowledge and possibilities, and far beyond the capabilities of the author.

\section{\sc Acknowledgement} I would like to express my thanks to Gary {\sc  Gundersen}, Katsuya {\sc Ishizaki}, and Kazuya {\sc Tohge} for valuable discussions on this topic several years ago, and once again to {Gary} for critically reading the first draft of this paper, thereby revealing several gaps and inaccuracies. My thanks also go to Katsutoshi {\sc Yamanoi}, whose reference to Corollary 4 in his paper \cite{Yamanoi} eventually allowed to conclude the proof of Theorem~\ref{AlgDependent}.

\end{document}